\pdfoutput=1
\RequirePackage{ifpdf}
\ifpdf % We are running pdfTeX in pdf mode
\documentclass[pdftex]{sigma}
\else
\documentclass{sigma}
\fi

\begin{document}

\allowdisplaybreaks

\renewcommand{\PaperNumber}{040}

\FirstPageHeading

\newcommand{\g}{\frak g}
\newcommand{\hg}{\hat {\frak g} }
\newcommand{\hn}{\hat {\frak n} }
\newcommand{\h}{\frak h}
\newcommand{\V}{\Cal V}
\newcommand{\hh}{\hat {\frak h} }
\newcommand{\n}{\frak n}
\newcommand{\Z}{\Bbb Z}
\newcommand{\N}{{\Bbb Z} _{> 0} }
\newcommand{\Zp} {\Z _ {\ge 0} }
\newcommand{\C}{\Bbb C}
\newcommand{\Q}{\Bbb Q}
\newcommand{\1}{\bf 1}
\newcommand{\wt}{{\rm {wt} }   }

\ShortArticleName{The Vertex Algebra $M(1)^+$ and Certain Af\/f\/ine Vertex Algebras of Level $-1$}

\ArticleName{The Vertex Algebra $\boldsymbol{M(1)^+}$\\ and Certain Af\/f\/ine Vertex Algebras of Level $\boldsymbol{-1}$}

\Author{Dra\v{z}en ADAMOVI\'{C} and Ozren PER\v{S}E}

\AuthorNameForHeading{D.~Adamovi\'{c} and O.~Per\v{s}e}

\Address{Faculty of Science, Department of
Mathematics, University of Zagreb,\\ Bijeni\v{c}ka cesta 30, 10000 Zagreb,
Croatia}

\Email{\href{mailto:adamovic@math.hr}{adamovic@math.hr}, \href{mailto:perse@math.hr}{perse@math.hr}}
\URLaddress{\url{http://web.math.pmf.unizg.hr/~adamovic/}}

\ArticleDates{Received March 09, 2012, in f\/inal form July 01, 2012; Published online July 08, 2012}

\Abstract{We give  a coset realization of the vertex operator
algebra $M(1)^+$ with central charge $\ell$. We realize $M(1) ^+$
as a commutant of certain af\/f\/ine vertex algebras of level $-1$ in
the vertex algebra $L_{C_{\ell} ^{(1)}}(-\tfrac{1}{2}\Lambda_0)
\otimes L_{C_{\ell} ^{(1)}}(-\tfrac{1}{2}\Lambda_0)$. We  show that
the simple vertex algebra $L_{C_{\ell} ^{(1)}}(-\Lambda_0)$ can be
(conformally) embedded into $L_{A_{2 \ell -1} ^{(1)}} (-\Lambda_0)$
and f\/ind the corresponding decomposition. We also study certain
coset subalgebras inside $L_{C_{\ell} ^{(1)}}(-\Lambda_0)$.}


\Classification{17B69; 17B67; 17B68; 81R10}

\section{Introduction}

In the last few years various types of $\mathcal{W}$-algebras have
been studied in the framework of vertex operator algebras (see
\cite{AdM-aim, AdM-cmp,Ar, DLY,KWak}).
In this paper we will be focused on ${\mathcal W}$-algebras which
admit coset realization. To any vertex algebra $V$ and its
subalgebra $U$, one can associate a new vertex algebra
\[ {\rm Com}(U,V)= \{ v \in V \; | \; u_n v =0 \mbox{ for all } u \in
U, \; n \geq 0 \}   \]
called the {\it commutant} (or {\it coset}) of $U$ in $V$. This is a
very important construction in the theory of vertex operator
algebras, because it gives a realization of a large family of
${\mathcal W}$-algebras.  Another important construction is the
{\it orbifold} construction, where a new vertex operator algebra is
obtained as invariants in a given vertex operator algebra with
respect to the f\/inite automorphism group. As we shall see in our
paper, some vertex algebras admit both realizations, coset and
orbifold.

In this paper we consider certain coset vertex algebras for vertex
algebras associated to af\/f\/ine Lie algebras. Let ${\frak g}$ be a
simple Lie algebra of type $X_{n}$, $\hat{\frak g}$ the associated
af\/f\/ine Lie algebra of type $X_{n}^{(1)}$, and $L_{X_{n}^{(1)}}(k
\Lambda _0)$ the simple vertex operator algebra associated to
$\hat{\frak g}$ of level $k \in \C$, $k \neq - h^{\vee}$. For $k,m
\in \N$, $L_{X_{n}^{(1)}}((k+m) \Lambda _0)$ is a subalgebra of
$L_{X_{n}^{(1)}}(k \Lambda _0) \otimes L_{X_{n}^{(1)}}(m \Lambda
_0)$, and one has the associated coset vertex operator algebra
\begin{gather} \label{coset-1}   {\rm Com}(L_{X_{n}^{(1)}}((k+m)
\Lambda _0) , L_{X_{n}^{(1)}}(k \Lambda _0) \otimes
L_{X_{n}^{(1)}}(m \Lambda _0)).
\end{gather}
Although there are no
precise general results (to the best of our knowledge) about the
structure of these cosets, it is believed that these vertex operator
algebras are f\/initely generated and rational.

In \cite{AP-2010},  we present a vertex-algebraic proof of the fact that  in  the case $k=1$ and af\/f\/ine Lie algebras of types $D_n ^{(1)}$ and $B_{n} ^{(1)}$,
  the coset (\ref{coset-1})  is isomorphic to the rational vertex operator algebra $V_L ^{+}$.

The situation is even more complicated for general $k,m \in \C$,
such that $k,m,k+m \neq - h^{\vee}$. Then one has the coset vertex
operator algebra
\[ {\rm Com}\big({\widetilde L}_{X_{n}^{(1)}}((k+m) \Lambda _0) ,
L_{X_{n}^{(1)}}(k \Lambda _0) \otimes L_{X_{n}^{(1)}}(m \Lambda
_0)\big), \]
where ${\widetilde L}_{X_{n}^{(1)}}((k+m) \Lambda _0)$ is a certain
af\/f\/ine vertex operator algebra associated to $X_{n}^{(1)}$ of level
$k+m$ (not necessarily simple). In this paper we identify some
special cases of such cosets, and it turns out that they are not
rational.

   The construction in \cite{AP-2010} is based on fermionic construction of vertex  algebras and certain conformal embeddings.
   In the present paper we use bosonic construction of vertex  algebras  and construct new conformal embeddings
   of af\/f\/ine vertex algebras at level $-1$. By applying the bosonic realization  of
the af\/f\/ine vertex  algebras $L_{A_1 ^{(1)}}(-\tfrac{1}{2}\Lambda_0)$
and $L_{C_{\ell} ^{(1)}}(-\tfrac{1}{2}\Lambda_0)$ (cf.~\cite{FF}) we
consider coset vertex   algebras
\begin{gather}   \mbox{Com}\big( L_{A_1 ^{(1)}}(- \Lambda_0),L_{A_1
^{(1)}}(-\tfrac{1}{2}\Lambda_0) \otimes L_{A_1
^{(1)}}(-\tfrac{1}{2}\Lambda_0)\big) \qquad \mbox{and} \nonumber\\
  {\rm Com} \big( {\widetilde L}_{C_{\ell} ^{(1)}}(-\Lambda_0),
L_{C_{\ell} ^{(1)}}(-\tfrac{1}{2}\Lambda_0) \otimes L_{C_{\ell}
^{(1)}}(-\tfrac{1}{2}\Lambda_0) \big). \label{coset-uvod}
\end{gather}
It is interesting that these cosets   have central charge $1$. We
show that these cosets are isomorphic to  $M(1) ^+$, where $M(1)$ is
the Heisenberg vertex operator algebra of rank $1$, and $M(1)^+$ is
the $\Z _2$-orbifold vertex algebra studied in~\cite{DN1}. The
structure theory of~$M(1) ^+$ shows that these cosets are irrational
vertex operator algebras and isomorphic to $W(2,4)$-algebra with
central charge $c=1$.

{\sloppy By combining results from \cite{A1} and the present paper, we
classify irreducible ordinary  $\widetilde{L}_{C_{\ell}
^{(1)}}(-\Lambda_0)$-modules. We believe that the (tensor) category
of  $\widetilde{L}_{C_{\ell} ^{(1)}}(-\Lambda_0)$-modules is
related  to the (tensor) category of $M(1) ^+$-modules. We plan to
address this correspondence in our forthcoming publications.

}

Generalizing (\ref{coset-uvod}), we use a natural realization of the
vertex operator algebra $L_{A_1 ^{(1)}}(- \Lambda_0) ^{\otimes
\ell}$ as a subalgebra of $L_{C_{\ell}
^{(1)}}(-\tfrac{1}{2}\Lambda_0) \otimes L_{C_{\ell}
^{(1)}}(-\tfrac{1}{2}\Lambda_0)$, and prove that
\[ \mbox{Com} \big(L_{A_1 ^{(1)}}(- \Lambda_0) ^{\otimes \ell},
L_{C_{\ell} ^{(1)}}(-\tfrac{1}{2}\Lambda_0) \otimes L_{C_{\ell}
^{(1)}}(-\tfrac{1}{2}\Lambda_0)\big) \]
is isomorphic to $M(1) ^+$, where $M(1)$ is the Heisenberg vertex
operator algebra of rank $\ell$.

Our construction is based on a new, interesting conformal  embedding
of af\/f\/ine vertex ope\-ra\-tor algebras at level~$-1$ which can be of
independent interest. We show that $L_{C_{\ell} ^{(1)}}(-\Lambda_0)$
is conformally embedded into $L_{A_{2 \ell -1} ^{(1)}} (-\Lambda_0)$
and that \[ L_{A_{2 \ell -1} ^{(1)}} (-\Lambda_0) \cong  L
_{C_{\ell} ^{(1)}} (-\Lambda_0) \oplus  L _{C_{\ell} ^{(1)}}
(-2\Lambda_0 + \Lambda_2),\]
which implies that $L _{C_{\ell} ^{(1)}} (-\Lambda_0)$ is a $\Z
_2$-orbifold of $L_{A_{2 \ell -1} ^{(1)}} (-\Lambda_0)$.

{\sloppy By using conformal embeddings we study certain categories of
$A_{2\ell-1} ^{(1)}$-modules of level~$-1$  from~\cite{AP} as
$C_{\ell} ^{(1)}$-modules. It turns out that  irreducible highest
weight $A_{2\ell-1} ^{(1)}$-modules $L_{ A_{2\ell -1} ^{(1)} } (-
(n+1) \Lambda_0 + n \Lambda_{1} )$ and $L_{ A_{2\ell -1} ^{(1)} } (-
(n+1) \Lambda_0 + n \Lambda_{2 \ell -1} )$ $(n \in {\N} )$ are also
irreducible as $C_{\ell} ^{(1)}$-modules. This result is an af\/f\/ine
analogue of the isomorphism of f\/inite-dimensional $
C_{\ell}$-modules:
\[ V_{A_{2\ell -1} } (n \omega_{2 \ell -1}) \cong V_{A_{2\ell -1} } (n
\omega_1) \cong V_{C_{ \ell  } } (n \omega_1).
\]

}

Using these conformal embeddings  we also  show that the coset
\[ {\rm Com} \big(   L _{A_1 ^{(1)}} (-\Lambda_0) ^{\otimes \ell}
,L_{C_{\ell} ^{(1)}} (-\Lambda_0) \big) \]
is isomorphic to $M(1) ^+$, where $M(1)$ is the Heisenberg vertex
operator algebra of rank $\ell -1$.

\section{Preliminaries}

Let $V$ be a vertex algebra \cite{B, FHL, FLM,LL}. For a subalgebra $U$ of $V$ denote by
\[
{\rm Com}(U,V)= \{ v \in V \; | \; u_n v =0 \mbox{ for all } u \in U,\;
n \geq 0 \} \]
 the commutant of $U$ in $V$ (cf.~\cite{FZ, GKO,LL}). Then, ${\rm Com}(U,V)$ is a
subalgebra of $V$ (also called coset vertex algebra).

Let ${\frak h}$ be a f\/inite-dimensional vector space equipped with a
nondegenerate symmetric bilinear form $\langle \cdot, \cdot
\rangle$, considered as an Abelian Lie algebra. Let $\hat{\frak h} =
{\frak h} \otimes {\C}[t,t^{-1}] \oplus {\C} K$ be its af\/f\/inization
with the center $K$. Then the free bosonic Fock space
$M(1) = S({\frak h} \otimes t^{-1} \C[t^{-1}])$ is a simple vertex
operator algebra of central charge $\ell = \dim {\frak h}$, with
Virasoro vector
\[
 \omega =\frac{1}{2} \sum _{i=1}^{\ell} h^{(i)}(-1)^2 {\1},\]
  where
$\{ h^{(1)}, \ldots , h^{(\ell)} \} $ is any orthonormal basis of
${\frak h}$ (cf.~\cite{FLM, LL}). We shall also use the
notation $M_{\frak h}(1)$ to emphasize the associated vector space
${\frak h}$.

Vertex algebra $M(1)$  has an order 2 automorphism which is lifted
from the map $h \mapsto -h$, for $h \in {\frak h}$. Denote by  $M(1)
^+$ (or $M_{\frak h}(1) ^+$) the subalgebra of invariants of that
automorphism. The irreducible modules for $M(1) ^+$ were classif\/ied
in \cite{DN1} and \cite{DN2}. For $\ell =1$, it was proved in
\cite{DG} that $M(1) ^+$ is generated by $\omega$ and one primary
vector of conformal weight $4$, so it is isomorphic to a
$W(2,4)$-algebra with central charge $1$.

Let $\frak g$ be the simple Lie algebra of type $X_{n}$, and
$\hat{\frak g}$ the associated af\/f\/ine Lie algebra of type
$X_{n}^{(1)}$. For any weight $\Lambda$ of $\hat{\frak g}$, denote
by $L_{X_{n} ^{(1)}}(\Lambda)$ the irreducible highest weight
$\hat{\frak g}$-module. Denote by~$\Lambda _i$, $i=0, \ldots, n$ the
fundamental weights of $\hat{\frak g}$ (cf.~\cite{K1}). We shall
also use the notation $V_{X_n}(\mu)$ for a~highest weight $\frak
g$-module of highest weight~$\mu$, and $\omega _i$, $i=1, \ldots, n$
for the fundamental weights of~$\frak g$.

{\sloppy For any $k \in \C$, denote by $N_{X_{n} ^{(1)}}(k \Lambda _0)$ the
generalized Verma $\hat{\frak g}$-module with highest weight~$k
\Lambda _0$. Then, $N_{X_{n} ^{(1)}}(k \Lambda _0)$ is a vertex
operator algebra of central charge $\frac{k \dim \frak
g}{k+h^{\vee}}$, for any $k \neq - h^{\vee}$, with Virasoro vector
obtained by Sugawara construction:
\begin{gather} \label{Vir-Sugawara}
\omega=\frac{1}{2(k+h^{\vee})}\sum_{i=1}^{\dim {\frak g}}
a^{i}(-1)b^{i}(-1){\bf 1},
\end{gather}
where $\{a^{i}\}_{i=1, \dots, \dim {\frak g}}$ is an arbitrary basis
of ${\frak g}$, and $\{b^{i}\}_{i=1, \dots, \dim {\frak g}}$  the
corresponding dual basis of~${\frak g}$ with respect to the
symmetric invariant bilinear form, normalized by the condition that
the length of the highest root is~$\sqrt{2}$ (cf.~\cite{FrB,FZ,K,LL,L}).

}

It follows that any quotient of $N_{X_{n} ^{(1)}}(k \Lambda _0)$ is
a vertex operator algebra, for $k \neq - h^{\vee}$. Specially,
$L_{X_{n} ^{(1)}}(k \Lambda _0)$ is a simple vertex operator
algebra, for any $k \neq - h^{\vee}$.

\section[Simple Lie algebras of type $C_{\ell}$ and $A_{2\ell -1}$]{Simple Lie algebras of type $\boldsymbol{C_{\ell}}$ and $\boldsymbol{A_{2\ell -1}}$}

Consider two $2 \ell$-dimensional vector spaces $A_1 = \oplus
_{i=1} ^{2 \ell}{\C} a_i ^+$, $ A_2 = \oplus _{i=1}^ {2 \ell} {\C}
a_i^-$ and let $A=A_1 {\oplus} A_2$. The Weyl algebra $W_{2 \ell}$ is
the complex associative algebra generated by $A$ with non-trivial
relations
\[
[a_i ^+,a_j^-]=\delta_{i,j} , \qquad 1 \le i,j \le 2 \ell.
\]
The normal ordering on $A$ is def\/ined by
\[
:\!xy\!:    = \frac 1 2 (xy+yx),  \qquad x,y \in A .
\]
Def\/ine
\[
e_{\epsilon_i-\epsilon_j}^A=  :\!a_i ^+ a_j^-\!:, \qquad
f_{\epsilon_i-\epsilon_j}^A=  :\!a_j ^+ a_i^-\!:, \qquad 1 \le i,j
\le 2 \ell, \quad i<j,
\]
and
\[
H_i=-:\!a_i ^+ a_i^-\!:,  \qquad 1 \le i \le 2 \ell.
\]
Then the Lie algebra $\frak g _1$ generated by the set
\[
\big\{ e_{\epsilon_i-\epsilon_j}^A, f_{\epsilon_i-\epsilon_j}^A \; \vert
\; 1 \le i,j \le 2 \ell, \ i < j \big\}
\]
is the simple Lie algebra of type $A_{2 \ell -1}$ (cf.~\cite{Bou}
and~\cite{FF}). The Cartan subalgebra $\frak h _1$ is spanned by
\[
\{ H_i - H_{i+1}  \; \vert \; 1 \le i \le 2 \ell -1 \}.
\]

Let $\theta$ be the automorphism of $W_{2 \ell}$ of order two given
by
\[  a_{i} ^+ \mapsto  a_{2\ell+1-i} ^{-}, \qquad a_{2\ell+1-i} ^-
\mapsto  a_{i} ^{+}, \qquad a_{i} ^-  \mapsto - a_{2\ell+1-i} ^+,
\qquad a_{2\ell+1-i} ^+ \mapsto - a_{i} ^-,\]
for $i =1, \dots, \ell$. Clearly, $\frak g _1$ is $\theta$-stable
and \[
\theta (e_{\epsilon_i-\epsilon_j}^A)=- e_{\epsilon_{2\ell
+1-j}-\epsilon_{2\ell +1-i}}^A, \qquad \theta
(e_{\epsilon_i-\epsilon_{2\ell +1-j}}^A)=
e_{\epsilon_{j}-\epsilon_{2\ell +1-i}}^A,
\]
and similarly for root vectors associated to negative roots.

The subalgebra $\frak g$ of $\frak g _1$ generated by
\begin{gather*}
  e_{2 \epsilon _{i}} =a_i ^+  a_{2\ell +1-i} ^-, \qquad
    f_{2 \epsilon _{i}}=a_i ^-  a_{2\ell +1-i} ^+,
  \\
  e_{\epsilon _{i}+ \epsilon _{j}}=\frac{1}{2}(a_i ^+
 a_{2\ell +1-j} ^- + a_j ^+  a_{2\ell +1-i} ^-), \qquad
 f_{\epsilon _{i}+ \epsilon _{j}}=\frac{1}{2}(a_i ^- a_{2\ell +1-j}
^+ + a_j ^-  a_{2\ell +1-i} ^+),   \\
  e_{\epsilon _{i}- \epsilon _{j}}=\frac{1}{2}(a_i ^+
 a_{j} ^- - a_{2\ell +1-j} ^+
 a_{2\ell +1-i} ^-), \qquad
f_{\epsilon _{i}- \epsilon _{j}}=\frac{1}{2}(a_j ^+
 a_{i} ^- - a_{2\ell +1-i} ^+
 a_{2\ell +1-j} ^-),
\end{gather*}
for $i,j =1, \dots, \ell$, $i<j$, is the simple Lie algebra of type
$C_{\ell}$. The Cartan subalgebra $\frak h$ is spanned by
\[
\{ H_i - H_{2\ell +1-i}  \; \vert \; 1 \le i \le  \ell \}.
\]
Clearly, $\theta$ acts as $1$ on $\frak g$. Furthermore, \[
e_{\epsilon _{1}+ \epsilon _{2}} ^*=\frac{1}{2}(a_1 ^+
 a_{2\ell -1} ^- - a_2 ^+  a_{2\ell} ^-) \in \frak g _1
\]
is a highest weight vector for $\frak g$, which generates the
irreducible $\frak g$-module $V_{C_{\ell}}(\omega _2)$. Clearly,
$\theta$~acts as $-1$ on $V_{C_{\ell}}(\omega _2)$. We obtain the
decomposition
\begin{gather} \label{decomp-fin1}
 \frak g _1 \cong \frak g \oplus V_{C_{\ell}}(\omega _2).
  \end{gather}
Since
\[ \dim V_{A_{2\ell -1}} (n \omega_1) = \dim V_{C_{\ell }} (n
\omega_1),\] one easily concludes that the irreducible $\frak g
_1$-module $V_{A_{2\ell -1}} (n \omega_1)$ remains irreducible when
restricted to $\frak g$. Thus,
\begin{gather} \label{decomp-fin2} V_{A_{2\ell -1} } (n \omega_1) \cong V_{C_{
\ell  } } (n \omega_1) \qquad \mbox{for} \quad n \in \Zp.
\end{gather}
Similarly
\begin{gather} \label{decomp-fin3} V_{A_{2\ell -1} } (n \omega_{2 \ell -1})
\cong V_{C_{ \ell  } } (n \omega_1) \qquad \mbox{for} \quad n \in \Zp.
\end{gather}
We will consider certain af\/f\/ine analogues of relations
(\ref{decomp-fin1})--(\ref{decomp-fin3}).

In what follows we shall need the following decompositions of $\frak
g$-modules:
\begin{gather}
V_{C_{\ell }} (\omega_2) \otimes V_{C_{\ell }} ( \omega_2) \cong
V_{C_{\ell }} (2 \omega_2)\oplus V_{C_{\ell }} ( \omega_1 +
\omega_3) \oplus V_{C_{\ell }} (\omega_4)
\nonumber \\
\hphantom{V_{C_{\ell }} (\omega_2) \otimes V_{C_{\ell }} ( \omega_2) \cong}{}
 \oplus V_{C_{\ell }} (2\omega_1) \oplus
V_{C_{\ell}} ( \omega_2) \oplus V_{C_{\ell }} (0)  \quad ( \ell \ge 4),  \nonumber \\
 V_{C_{3}} (\omega_2) \otimes V_{C_{3}} ( \omega_2) \cong
 V_{C_{3}} (2 \omega_2)\oplus V_{C_{3}} ( \omega_1 +
\omega_3) \oplus V_{C_{3}} (2\omega_1) \oplus
V_{C_{3}} ( \omega_2) \oplus V_{C_{3}} (0), \nonumber \\
 V_{C_{2}}(\omega _2) \otimes V_{C_{2}}(\omega _2) \cong
V_{C_{2}}(2 \omega _2) \oplus V_{C_{2}}(2 \omega _1) \oplus V_{C_{2}}(0). \label{tens-pr-decomp}
 \end{gather}

\section{Weyl vertex algebras and symplectic af\/f\/ine Lie algebras}

The Weyl algebra $W_{\ell}( \tfrac{1}{2} + {\Z})$ is a complex
associative algebra generated by
\[  a ^{\pm} _{i}(r),  \qquad   r  \in \tfrac{1}{2} + {\Z}, \quad 1 \le i
\le \ell
\] with non-trivial relations
\begin{gather*}
 [a_{i} ^{+}(r)  , a_{j} ^{-}(s)  ]  =\delta_{r+s,0} \delta_{i,j},
\end{gather*}
where $r, s \in {\tfrac{1}{2}}+ {\Z}$, $ i, j \in \{1, \dots, \ell
\}$.

 Let $M_{\ell}$ be the irreducible $W_{\ell}( \tfrac{1}{2} + {\Z})$-module generated by
 the
cyclic vector ${\1}$ such that
\[ a^{\pm} _i  (r) {\1}     =  0 \qquad \mbox{for} \quad r > 0 ,  \quad 1
\le i \le \ell.
\]

Def\/ine the following   f\/ields on $M_{\ell}$
\[   a_i ^{\pm }(z) = \sum_{ n \in   {\Z}
 } a_i ^{\pm }(n+ \tfrac{1}{2} )  z ^{-n- 1}.
 \]
 The f\/ields $a_{i} ^{\pm }(z)$, $i =1, \dots, \ell$ generate on $M_{\ell}$  the
unique structure of a simple vertex algebra (cf.~\cite{FrB, K}). Let us denote the corresponding vertex operator by $Y$.

We have the following Virasoro vector in $M_{\ell}$:
\begin{gather} \label{vir-bos}  \omega=\frac{1}{2} \sum_{i=1} ^{\ell} \big( a_i
^-  (- \tfrac{3}{2} ) a_i ^{+}  (-\tfrac{1}{2}  ) -a_i ^+  (-  \tfrac{3}{2} )
a_i ^{-} (-\tfrac{1}{2}  ) \big) {\1}.
\end{gather}
Let $Y(\omega,z) = \sum\limits_{n \in {\Z} } L(n) z ^{-n-2}$.
Then $M_{\ell}$ is $\tfrac{1}{2}{\Zp}$-graded with respect to
$L(0)$:
\[ M_{\ell} := \bigoplus_{ m \in \tfrac{1}{2} \Zp} M_{\ell} (m),
\qquad M_{\ell} (m)= \{ v \in M_{\ell} \; \vert \;  L(0) v = m v \}.\]
Note that $M_{\ell} (0) = {\C} {\1}$. For $ v \in M_{\ell} (m)$ we
shall write $\wt (v) = m$.

The following result is well-known.

\begin{theorem}[\cite{FF}] \label{ff}
We have
\[  M_{\ell} \cong L_{C_{\ell} ^{(1)}}( -\tfrac{1}{2}\Lambda _{0})
\bigoplus L_{C_{\ell} ^{(1)}}( -\tfrac{3}{2}\Lambda _{0} +
\Lambda_1). \]
In the case $\ell =1$ we have \[  M_1 \cong L_{A_1 ^{(1)}}(
-\tfrac{1}{2}\Lambda _{0}) \bigoplus L_{A_1 ^{(1)}}(
-\tfrac{3}{2}\Lambda _{0} + \Lambda_1). \]
\end{theorem}

\begin{remark} The highest weights of modules from Theorem \ref{ff}
are admissible in the sense of~\cite{KW}. Representations of vertex
operator algebras associated to af\/f\/ine Lie algebras of type~$A_1
^{(1)}$ and~$C_{\ell} ^{(1)}$ with admissible highest weights were
studied in~\cite{AM} and~\cite{A-1994}.
\end{remark}

We shall now consider the vertex algebra $M_{2 \ell}$ and its
subalgebra $L_{C_{\ell} ^{(1)}}( -\tfrac{1}{2}\Lambda _{0}) \otimes
L_{C_{\ell} ^{(1)}}( -\tfrac{1}{2}\Lambda _{0})$. For $\ell =1$,
$L_{A_1 ^{(1)}}(-\tfrac{1}{2}\Lambda_0) \otimes L_{A_1
^{(1)}}(-\tfrac{1}{2}\Lambda_0)$ is a subalgebra of $M_{2}$.

Let $\theta : M_{2 \ell} \rightarrow M_{2 \ell}$ be the automorphism
of order two of  the vertex algebra $M_{2 \ell}$ which is lifted
from the following automorphism of the Weyl algebra
\begin{alignat*}{3}
& a_{i} ^+ (s) \mapsto  a_{2\ell+1-i} ^{-} (s), \qquad
&& a_{2\ell+1-i} ^- (s) \mapsto  a_{i} ^{+} (s), &  \\
& a_{i} ^-(s) \mapsto - a_{2\ell+1-i} ^+ (s), \qquad && a_{2\ell+1-i}
^+(s) \mapsto - a_{i} ^- (s), &
 \end{alignat*}
for $i =1, \dots, \ell$ and  $s \in \tfrac{1}{2} + {\Z}$.

If we have a subalgebra $U \subset M_{2 \ell}$ which is
$\theta$-stable, we def\/ine
\[
U ^0 = \{ u \in U \; \vert \; \theta (u) = u \}, \qquad U ^1 = \{ u
\in U \; \vert \; \theta (u) = -u \}.
\]

Def\/ine the following vectors in $M_{2 \ell}$:
\begin{alignat*}{3}
&b_{i} ^+ = \frac{1}{\sqrt{2}}( a_{i} ^+ (-\tfrac{1}{2}) +
a_{2\ell+1-i} ^- (-\tfrac{1}{2}) ) {\1},\qquad
&& b_{2\ell+1-i} ^+ = \frac{\sqrt{-1}}{\sqrt{2}}( a_{i} ^+ (-\tfrac{1}{2}) - a_{2\ell+1-i} ^- (-\tfrac{1}{2}) ) {\1}, &  \\
& b_{i} ^- = \frac{1}{\sqrt{2}}( a_{i} ^- (-\tfrac{1}{2}) -
a_{2\ell+1-i}^+ (-\tfrac{1}{2}) ) {\1}, \qquad &&b_{2\ell+1-i} ^- =
\frac{-\sqrt{-1}}{\sqrt{2}}( a_{i} ^- (-\tfrac{1}{2}) +
a_{2\ell+1-i} ^+ (-\tfrac{1}{2}) ) {\1}, &
\end{alignat*} for $i =1, \dots, \ell$. Then the subalgebra generated by $b_i
^+$, $b _i ^-$ (resp.\ $b_{2\ell+1-i} ^+$, $b _{2\ell+1-i} ^-$), for $i
=1, \dots, \ell$, is isomorphic to $M_{\ell}$.
Since
\[ \theta (b_i ^{\pm}) = b_i ^{\pm}, \qquad  \theta (b_{2\ell+1-i}
^{\pm}) = - b_{2\ell+1-i} ^{\pm},\]
for $i =1, \dots, \ell$, we have \[ M_{2 \ell}^{0} \cong M_{\ell}
\otimes M_{\ell} ^0 \cong M_{\ell} \otimes L_{C_{\ell} ^{(1)}}(
-\tfrac{1}{2}\Lambda _{0}),\] and for $\ell =1$: \[ M_{2}^{0} \cong
M_{1} \otimes M_{1} ^0 \cong M_{1} \otimes L_{A_1 ^{(1)}}(
-\tfrac{1}{2}\Lambda _{0}).\]

\section[Commutant of $L_{A_1 ^{(1)}}(- \Lambda_0) ^{\otimes \ell}$ in  $L_{C_{\ell} ^{(1)}}(-\frac{1}{2}\Lambda_0)
\otimes L_{C_{\ell} ^{(1)}}(-\frac{1}{2}\Lambda_0)$]{Commutant of $\boldsymbol{L_{A_1 ^{(1)}}(- \Lambda_0) ^{\otimes \ell}}$ in  $\boldsymbol{L_{C_{\ell} ^{(1)}}(-\frac{1}{2}\Lambda_0)
\otimes L_{C_{\ell} ^{(1)}}(-\frac{1}{2}\Lambda_0)}$}

In this section we use Weyl vertex algebra $M_{2 \ell}$ to study
certain subalgebras of $L_{C_{\ell} ^{(1)}}(-\tfrac{1}{2}\Lambda_0)
\otimes L_{C_{\ell} ^{(1)}}(-\tfrac{1}{2}\Lambda_0)$. We determine
the commutants
\[
{\rm Com} \big( L_{A_1 ^{(1)}}(- \Lambda_0) ^{\otimes
\ell},L_{C_{\ell} ^{(1)}}(-\tfrac{1}{2}\Lambda_0) \otimes
L_{C_{\ell} ^{(1)}}(-\tfrac{1}{2}\Lambda_0)\big)
\]
and \[
{\rm Com} \big( L_{A_1 ^{(1)}}(-\Lambda_0),L_{A_1
^{(1)}}(-\tfrac{1}{2}\Lambda_0) \otimes L_{A_1
^{(1)}}(-\tfrac{1}{2}\Lambda_0) \big).
\]

Let $U$ be the vertex subalgebra  of $M_{\ell} ^0 \otimes M_{\ell}
^0 \subset M_{2 \ell} ^0$ generated by
\begin{gather}
e^{(i)}=a_i ^+ (-\tfrac{1}{2}) a_{
2\ell +1-i} ^-( -\tfrac{1}{2}) {\1}, \qquad  f^{(i)}=a_ i ^-
(-\tfrac{1}{2}) a_{2\ell +1-i} ^+(
-\tfrac{1}{2}) {\1} \qquad \mbox{and} \nonumber \\
 h^{(i)}= (-a_ i ^+ (-\tfrac{1}{2}) a_i^-( -\tfrac{1}{2}) +
a_{2\ell +1-i} ^+ (-\tfrac{1}{2}) a_{2\ell +1-i} ^-( -\tfrac{1}{2})){\1}, \label{gen-1-higher}
 \end{gather}
$i =1, \dots, \ell$. It is clear that $U$ is isomorphic to the tensor product
\[
 \underbrace{ L_{A_1 ^{(1)}}( -\Lambda _{0}) \otimes \cdots \otimes  L_{A_1 ^{(1)}}( -\Lambda _{0})}_{\ell \ \  \mbox{times} }
\]
of $\ell$ copies of the af\/f\/ine vertex algebra $ L_{A_1 ^{(1)}}(-\Lambda _{0})$.

Def\/ine also
 \[
 H ^{(i)} =   (a_i ^+ (-\tfrac{1}{2}) a_i ^- (-\tfrac{1}{2})+a_{2\ell +1-i} ^+ (-\tfrac{1}{2})
 a_{2\ell +1-i} ^- (-\tfrac{1}{2}) ) {\1}, \qquad i=1, \dots,
 \ell.
 \]
 Then
 \[ H ^{(i)} \in \mbox{Com} ( U, M_{2\ell}), \qquad i=1, \dots,  \ell.
 \]

 Let \[ {\frak h}= \bigoplus_{i =1 } ^{ \ell} {\C} H ^{(i)}.\] Then ${\frak h}$ can be
 considered as an Abelian Lie algebra, and the components of the vertex operators
 \[
 Y(h,z) = \sum_{\n \in {\Z} } h(n) z ^{-n-1}, \qquad h \in {\frak  h},
 \]
 def\/ine a representation of the associated Heisenberg algebra $\widehat{\frak h}$.
 Moreover, ${\frak h}$ generates the subalgebra of the~$M_{2 \ell}$ which is isomorphic to the
 Heisenberg vertex algebra $M_{\frak h} (1)$ with central charge~$\ell$.
 We also note that $\langle H ^{(i)}, H ^{(j)} \rangle = H ^{(i)}(1)H ^{(j)}= -2 \delta _{ij}$.

Using relations (\ref{gen-1-higher}) one can show that the Virasoro
vector (\ref{vir-bos}) (in $M_{2 \ell}$) can be written in a~form:
\[ \omega =  \omega_1   + \omega_2 \]
where \[\omega_1 =\frac{1}{2} \sum_{i = 1 } ^{\ell} \big(e
^{(i)}(-1)f^{(i)}(-1)+f ^{(i)} (-1)e^{(i)}(-1)+ \tfrac{1}{2} h ^{(i)}
(-1)^2\big){\1}
\] is the  Virasoro vector  in $U$  and
\[\omega_2 = -\frac{1}{4}\sum_{i = 1 } ^{\ell} H ^{(i) }(-1) ^{2}
{\1}\] is the Virasoro vector   in $M_{\frak h} (1)$.
For $i=1,2$ let
$Y(\omega_i , z) = \sum\limits_{n\in {\Z} } L_i (n) z ^{-n-2}$.

\begin{proposition} \label{com-bos-higher}
We have:
\[
{\rm Com} (U, M_{2\ell}) = M_{\frak h} (1).
\]
\end{proposition}

\begin{proof}
It is clear that
 $W= \mbox{Com} ( U, M_{2\ell})$ contains a subalgebra isomorphic to $M_{\frak h}
(1)$.

Assume now that $M_{\frak h} (1) \ne W$. Then there is a vector
$ w \in W$, ${\wt} (w) > 0$ such that
\[ H ^{(i)}(n) w =  \delta_{n,0} \lambda _{(i)} w, \qquad n \in {\Zp},
\quad \lambda _{(i)} \in {\Z}, \quad i = 1, \dots, \ell.
\]
(Note that each  $H^{(i)}(0)$ acts semisimply on $M_{2\ell}$ with
eigenvalues in ${\Z}$.)
By the def\/inition of $W$ we have that $L_1 (0) w = 0$. Therefore
\[ L(0) w = L_2(0) w = - \frac{1}{4} \sum_{i = 1 } ^{\ell} (\lambda
_{(i)})  ^2 w. \] This contradicts the fact that $\wt (w) > 0$. So,
$W= M_{\frak h} (1)$.
\end{proof}

Since $\theta (H^{(i)}) = -H^{(i)}$ for $i=1, \dots, \ell$ we have
that $M_{\frak h} (1)^+ =M_{\frak h}(1) ^0$ is a subalgebra of
$M_{\ell} ^0 \otimes M_{\ell} ^0$.
 Therefore the  vertex algebra $L_{C_{\ell} ^{(1)}}(-\tfrac{1}{2}\Lambda_0) \otimes L_{C_{\ell} ^{(1)}}(-\tfrac{1}{2}\Lambda_0)$ contains a subalgebra isomorphic
to $U \otimes M_{\frak h} (1) ^+$. By using Proposition~\ref{com-bos-higher} we obtain the following theorem.

\begin{theorem} \label{coset-bos-r-higher}We have:
\[ {\rm Com} \big( L_{A_1 ^{(1)}}(- \Lambda_0) ^{\otimes
\ell},L_{C_{\ell} ^{(1)}}(-\tfrac{1}{2}\Lambda_0) \otimes
L_{C_{\ell} ^{(1)}}(-\tfrac{1}{2}\Lambda_0) \big) \cong M_{\frak h} (1)
^+.\]
\end{theorem}

In the case $\ell =1$, ${\frak h}$ is a one-dimensional vector space
${\frak h}= {\C} H$, where
\[ H = \big(a_1 ^+ (-\tfrac{1}{2}) a_1 ^-
(-\tfrac{1}{2})+a_2 ^+ (-\tfrac{1}{2}) a_2 ^- (-\tfrac{1}{2}) \big)
{\1}.\]

\begin{theorem} \label{coset-bos-r} We have:
\[ {\rm Com} \big( L_{A_1 ^{(1)}}(-\Lambda_0),L_{A_1
^{(1)}}(-\tfrac{1}{2}\Lambda_0) \otimes L_{A_1
^{(1)}}(-\tfrac{1}{2}\Lambda_0) \big) \cong M_{\frak h} (1) ^+.
\]
\end{theorem}

In the next section we generalize Theorem \ref{coset-bos-r} in
another way.

\section[Commutant of level $-1$ type $C_{\ell} ^{(1)}$  affine vertex algebra in  $L_{C_{\ell} ^{(1)}}(-\frac{1}{2}\Lambda_0)
\otimes L_{C_{\ell} ^{(1)}}(-\frac{1}{2}\Lambda_0)$]{Commutant of level $\boldsymbol{-1}$ type $\boldsymbol{C_{\ell} ^{(1)}}$
af\/f\/ine vertex algebra\\ in  $\boldsymbol{L_{C_{\ell} ^{(1)}}(-\frac{1}{2}\Lambda_0)
\otimes L_{C_{\ell} ^{(1)}}(-\frac{1}{2}\Lambda_0)}$}

In this section we study another subalgebra of $L_{C_{\ell}
^{(1)}}(-\frac{1}{2}\Lambda_0) \otimes L_{C_{\ell}
^{(1)}}(-\frac{1}{2}\Lambda_0)$. The subalgebra of $M_{\ell} ^0
\otimes M_{\ell} ^0 \subset M_{2 \ell} ^0$ generated by
\begin{gather}
  e_{2 \epsilon _{i}} (=e^{(i)})=a_i ^+ (-\tfrac{1}{2}) a_{2\ell +1-i} ^-( -\tfrac{1}{2})
 {\1}, \qquad   f_{2 \epsilon _{i}}(=f^{(i)})=a_i ^- (-\tfrac{1}{2}) a_{2\ell +1-i} ^+( -\tfrac{1}{2}) {\1},
 \nonumber \\
  e_{\epsilon _{i}+ \epsilon _{j}}=\tfrac{1}{2}\big(a_i ^+
(-\tfrac{1}{2}) a_{2\ell +1-j} ^-( -\tfrac{1}{2}) + a_j ^+
(-\tfrac{1}{2}) a_{2\ell +1-i} ^-( -\tfrac{1}{2})\big){\1}, \nonumber \\
  f_{\epsilon _{i}+ \epsilon _{j}}=\tfrac{1}{2}\big(a_i ^-
(-\tfrac{1}{2}) a_{2\ell +1-j} ^+( -\tfrac{1}{2}) + a_j ^-
(-\tfrac{1}{2}) a_{2\ell +1-i} ^+( -\tfrac{1}{2})\big){\1}, \nonumber \\
  e_{\epsilon _{i}- \epsilon _{j}}=\tfrac{1}{2}\big(a_i ^+
(-\tfrac{1}{2}) a_{j} ^-( -\tfrac{1}{2}) - a_{2\ell +1-j} ^+
(-\tfrac{1}{2}) a_{2\ell +1-i} ^-( -\tfrac{1}{2})\big){\1}, \nonumber \\
  f_{\epsilon _{i}- \epsilon _{j}}=\tfrac{1}{2}\big(a_j ^+
(-\tfrac{1}{2}) a_{i} ^-( -\tfrac{1}{2}) - a_{2\ell +1-i} ^+
(-\tfrac{1}{2}) a_{2\ell +1-j} ^-( -\tfrac{1}{2})\big){\1}, \nonumber
\\
  \hphantom{f_{\epsilon _{i}- \epsilon _{j}}=}{}
    \ \mbox{for}\quad i,j =1, \dots, \ell, \ \ i<j,\label{gen-1-C}
  \end{gather}
is a level $-1$ af\/f\/ine vertex operator algebra associated to
$C_{\ell} ^{(1)}$. We denote it by ${\widetilde L}_{C_{\ell}
^{(1)}}(- \Lambda_0)$.

Let
\[ H = \sum_{i = 1 } ^{\ell} H ^{(i)} =\sum _{i=1}^{2 \ell} a_i ^+
(-\tfrac{1}{2}) a_i ^- (-\tfrac{1}{2}){\1}. \]

Set ${\h}_1 = {\C} H \subset {\h}$. Let $M_{ {\h}_1} (1)$ be the
Heisenberg vertex algebra generated by $H$. Clearly $\langle H, H
\rangle = -2 \ell$.

Since  $H \in \mbox{Com} ({\widetilde L}_{C_{\ell} ^{(1)}}(-
\Lambda_0), M_{2 \ell})$, we have that
  $M_{2 \ell}$ contains a
subalgebra isomorphic to $ {\widetilde L}_{C_{\ell} ^{(1)}}(-
\Lambda_0) \otimes M_{{\h}_1 }(1)$.

\begin{proposition} \label{Vir-decomp-C}
The Virasoro vector \eqref{vir-bos} in $M_{2 \ell}$ can be written
in a form:
\[ \omega = \omega_1 + \omega_2,
\]
where $ \omega_1$ is the  Virasoro vector \eqref{Vir-Sugawara} in
${\widetilde L}_{C_{\ell} ^{(1)}}(- \Lambda_0)$ obtained by the
Sugawara construction and
\[ \omega_2 = -\frac{1}{4 \ell} H(-1) ^{2} {\1}\] is the Virasoro
vector in $M_{ {\h}_1} (1)$.
\end{proposition}

\begin{proof} Formula (\ref{Vir-Sugawara}) implies that
\begin{gather} \omega_1 = \frac{1}{2 \ell} \Bigg(
\sum_{i=1}^{\ell} (e_{2 \epsilon _{i}}(-1)f_{2 \epsilon
_{i}}(-1)+f_{2 \epsilon
_{i}}(-1)e_{2 \epsilon _{i}}(-1)){\1} \nonumber \\
\hphantom{\omega_1 =}{}  +2 \sum_{i,j=1 \atop i<j}^{\ell} (e_{\epsilon _{i}+ \epsilon
_{j}}(-1)f_{\epsilon _{i}+ \epsilon _{j}}(-1)+ f_{\epsilon _{i}+
\epsilon _{j}}(-1)e_{\epsilon _{i}+ \epsilon _{j}}(-1)){\1} \nonumber \\
\hphantom{\omega_1 =}{}
  +2 \sum_{i,j=1 \atop i<j}^{\ell} (e_{\epsilon _{i}- \epsilon
_{j}}(-1)f_{\epsilon _{i}- \epsilon _{j}}(-1)+ f_{\epsilon _{i}-
\epsilon _{j}}(-1)e_{\epsilon _{i}- \epsilon _{j}}(-1)){\1}
  +\frac{1}{2}\sum_{i=1}^{\ell} h_{2 \epsilon _{i}}(-1)^2 {\1}
\Bigg),\label{proof-Vir-1}
 \end{gather}
where
\begin{gather*} h_{2 \epsilon _{i}} (=h^{(i)})= \big({-}a_i ^+ (-\tfrac{1}{2}) a_i^-(
-\tfrac{1}{2})+ a_{2\ell +1-i} ^+ (-\tfrac{1}{2}) a_{2\ell +1-i}^-(
-\tfrac{1}{2}) \big){\1}.
 \end{gather*}
It follows from relations (\ref{gen-1-C}) that
\begin{gather}
 (e_{2 \epsilon _{i}}(-1)f_{2 \epsilon
_{i}}(-1)+f_{2
\epsilon _{i}}(-1)e_{2 \epsilon _{i}}(-1)){\1}\nonumber\\
\qquad{}
 = 2 a_i ^+ (-\tfrac{1}{2})a_i ^- (-\tfrac{1}{2}) a_{2\ell +1-i}
^+ (-\tfrac{1}{2}) a_{2\ell +1-i} ^- (-\tfrac{1}{2}){\1} \nonumber \\
\quad\qquad{} + a_i ^- (-\tfrac{3}{2})a_i ^+ (-\tfrac{1}{2}){\1}  + a_{2\ell
+1-i} ^- (-\tfrac{3}{2})a_{2\ell +1-i} ^+ (-\tfrac{1}{2}){\1} \nonumber \\
\quad\qquad{} -a_i ^+ (-\tfrac{3}{2})a_i ^- (-\tfrac{1}{2}){\1}  - a_{2\ell
+1-i} ^+ (-\tfrac{3}{2})a_{2\ell +1-i} ^- (-\tfrac{1}{2}){\1},\label{proof-Vir-2}
\\
 2(e_{\epsilon _{i}+ \epsilon
_{j}}(-1)f_{\epsilon _{i}+ \epsilon _{j}}(-1)+ f_{\epsilon _{i}+
\epsilon _{j}}(-1)e_{\epsilon _{i}+
\epsilon _{j}}(-1)){\1}  \nonumber\\
\qquad{}
 = a_i ^+ (-\tfrac{1}{2})a_i ^- (-\tfrac{1}{2}) a_{2\ell +1-j} ^+
(-\tfrac{1}{2}) a_{2\ell +1-j} ^- (-\tfrac{1}{2}){\1} \nonumber \\
\quad\qquad{} + a_i ^+ (-\tfrac{1}{2})a_j ^- (-\tfrac{1}{2}) a_{2\ell +1-j} ^-
(-\tfrac{1}{2})
a_{2\ell +1-i} ^+ (-\tfrac{1}{2}){\1} \nonumber \\
\quad\qquad{} + a_i ^- (-\tfrac{1}{2})a_j ^+ (-\tfrac{1}{2}) a_{2\ell +1-i} ^-
(-\tfrac{1}{2}) a_{2\ell +1-j} ^+ (-\tfrac{1}{2}){\1} \nonumber \\
\quad\qquad{} + a_j ^+ (-\tfrac{1}{2})a_j ^- (-\tfrac{1}{2}) a_{2\ell +1-i} ^+
(-\tfrac{1}{2})
a_{2\ell +1-i} ^- (-\tfrac{1}{2}){\1} \nonumber \\
\quad \qquad{} +\frac{1}{2} \Big( a_i ^- (-\tfrac{3}{2}) a_i ^+
(-\tfrac{1}{2}){\1} + a_j ^- (-\tfrac{3}{2}) a_j ^+
(-\tfrac{1}{2}){\1} + a_{2\ell +1-i} ^- (-\tfrac{3}{2}) a_{2\ell
+1-i} ^+
(-\tfrac{1}{2}){\1}  \nonumber \\
\quad \qquad{} +a_{2\ell +1-j} ^- (-\tfrac{3}{2}) a_{2\ell +1-j} ^+
(-\tfrac{1}{2}){\1} - a_i ^+ (-\tfrac{3}{2}) a_i ^-
(-\tfrac{1}{2}){\1} - a_j ^+ (-\tfrac{3}{2}) a_j ^-
(-\tfrac{1}{2}){\1} \nonumber \\
\quad\qquad{} - a_{2\ell +1-i} ^+ (-\tfrac{3}{2}) a_{2\ell
+1-i}^- (-\tfrac{1}{2}){\1} - a_{2\ell +1-j} ^+ (-\tfrac{3}{2})
a_{2\ell
+1-j} ^- (-\tfrac{1}{2}){\1} \Big),\label{proof-Vir-3}
  \end{gather}
and \begin{gather}
 2(e_{\epsilon _{i}- \epsilon
_{j}}(-1)f_{\epsilon _{i}- \epsilon _{j}}(-1)+ f_{\epsilon _{i}-
\epsilon _{j}}(-1)e_{\epsilon _{i}-
\epsilon _{j}}(-1)){\1} \nonumber \\
\qquad{} = a_i ^+ (-\tfrac{1}{2})a_i ^- (-\tfrac{1}{2}) a_{j} ^+
(-\tfrac{1}{2}) a_{j} ^- (-\tfrac{1}{2}){\1} \nonumber \\
\quad\qquad{} - a_j ^+ (-\tfrac{1}{2})a_i ^- (-\tfrac{1}{2}) a_{2\ell +1-j} ^+
(-\tfrac{1}{2})
a_{2\ell +1-i} ^- (-\tfrac{1}{2}){\1} \nonumber \\
\quad\qquad{} - a_i ^+ (-\tfrac{1}{2})a_j ^- (-\tfrac{1}{2}) a_{2\ell +1-i} ^+
(-\tfrac{1}{2}) a_{2\ell +1-j} ^- (-\tfrac{1}{2}){\1}  \nonumber \\
\quad\qquad{} +a_{2\ell +1-i} ^+ (-\tfrac{1}{2})a_{2\ell +1-i} ^-
(-\tfrac{1}{2}) a_{2\ell +1-j} ^+ (-\tfrac{1}{2})
a_{2\ell +1-j} ^- (-\tfrac{1}{2}){\1} \nonumber \\
\quad\qquad{} +\frac{1}{2} \Big( a_i ^- (-\tfrac{3}{2}) a_i ^+
(-\tfrac{1}{2}){\1} + a_j ^- (-\tfrac{3}{2}) a_j ^+
(-\tfrac{1}{2}){\1} + a_{2\ell +1-i} ^- (-\tfrac{3}{2}) a_{2\ell
+1-i} ^+
(-\tfrac{1}{2}){\1}  \nonumber \\
\quad\qquad{} +a_{2\ell +1-j} ^- (-\tfrac{3}{2}) a_{2\ell +1-j} ^+
(-\tfrac{1}{2}){\1} - a_i ^+ (-\tfrac{3}{2}) a_i ^-
(-\tfrac{1}{2}){\1} - a_j ^+ (-\tfrac{3}{2}) a_j ^-
(-\tfrac{1}{2}){\1} \nonumber \\
\quad\qquad{} - a_{2\ell +1-i} ^+ (-\tfrac{3}{2}) a_{2\ell
+1-i}^- (-\tfrac{1}{2}){\1} - a_{2\ell +1-j} ^+ (-\tfrac{3}{2})
a_{2\ell
+1-j} ^- (-\tfrac{1}{2}){\1} \Big), \label{proof-Vir-4}
  \end{gather}
for all $i,j =1, \dots, \ell$, $i<j$. Using relations
(\ref{proof-Vir-1}), (\ref{proof-Vir-2}), (\ref{proof-Vir-3}) and
(\ref{proof-Vir-4}) one can obtain
\begin{gather*} \omega _{1}=
\frac{1}{2} \sum_{i=1} ^{2 \ell} \big( a_i ^- (- \tfrac{3}{2}) a_i ^{+}
(-\tfrac{1}{2} ) -a_i ^+ (- \tfrac{3}{2}) a_i ^{-} (-\tfrac{1}{2} ) \big)
{\1} + \frac{1}{4 \ell} H(-1) ^{2} {\1},
 \end{gather*}
which implies the claim of proposition.
\end{proof}

Using Proposition \ref{Vir-decomp-C} and applying similar arguments
as in the proof of Proposition \ref{com-bos-higher}, we obtain:

\begin{proposition} \label{com-bos-C} We have:
\[ {\rm Com} \big({\widetilde L}_{C_{\ell} ^{(1)}}(- \Lambda_0), M_{2
\ell}\big) = M_ { {\h}_1} (1).\]
\end{proposition}

 Since $\theta (H) = -H$, we have that $M_ { {\h}_1} (1)^+
=M_ { {\h}_1} (1) ^0$ is a subalgebra of $M_{\ell} ^0 \otimes
M_{\ell} ^0$.
 Therefore the  vertex algebra $L_{C_{\ell} ^{(1)}}(-\tfrac{1}{2}\Lambda_0) \otimes
 L_{C_{\ell} ^{(1)}}(-\tfrac{1}{2}\Lambda_0)$ contains a subalgebra isomorphic
to ${\widetilde L}_{C_{\ell} ^{(1)}}(- \Lambda_0) \otimes M_ {
{\h}_1} (1) ^+$. By using Proposition \ref{com-bos-C} we obtain the
following theorem.

\begin{theorem} \label{coset-bos-r-C}We have:
\[ {\rm Com} \big( {\widetilde L}_{C_{\ell} ^{(1)}}(-\Lambda_0),
L_{C_{\ell} ^{(1)}}(-\tfrac{1}{2}\Lambda_0) \otimes L_{C_{\ell}
^{(1)}}(-\tfrac{1}{2}\Lambda_0) \big) \cong M_ { {\h}_1}  (1) ^+.\]
\end{theorem}

\section[The classification of ordinary modules for  ${\widetilde L}_{C_{\ell} ^{(1)}}(-
\Lambda_0)$]{The classif\/ication of ordinary modules for  $\boldsymbol{{\widetilde L}_{C_{\ell} ^{(1)}}(-
\Lambda_0)}$} \label{sect.class.C}

In this section we obtain a classif\/ication of irreducible
${\widetilde L}_{C_{\ell} ^{(1)}}(- \Lambda_0)$-modules, which we
use in the following sections.

The following vertex algebra was considered in \cite{A1}: Let $\ell
\ge 3$ and \[
V_{C_{\ell} ^{(1)}} (-\Lambda_0)= \frac{N_{C_{\ell} ^{(1)}} (-\Lambda_0)}{\langle \Delta_3 (-1) {\1} \rangle },
\]
where $\langle \Delta_3 (-1) {\1} \rangle$ is the ideal generated by the
singular vector $\Delta_3 (-1) {\1}$, such that $\Delta_3(-1)$ is
given by the following determinant:
\begin{gather*}
 \Delta_3(-1) =
 \left| \begin{array}{ccc} e_{ 2\epsilon_1} (-1) & e_{
\epsilon_1 + \epsilon_2} (-1) &  e_{ \epsilon_1 + \epsilon_3} (-1)
\\ e_{ \epsilon_1 + \epsilon_2} (-1) & e_{ 2
\epsilon_2} (-1)&  e_{ \epsilon_2 + \epsilon_3} (-1)
\\e_{ \epsilon_1
+ \epsilon_3} (-1) & e_{ \epsilon_2 + \epsilon_3} (-1) & e_{ 2
\epsilon_3} (-1)
\end{array} \right|.
 \end{gather*}
Using relations (\ref{gen-1-C}), one can easily check that
\begin{gather*}
 \left| \begin{array}{ccc} e_{ 2\epsilon_1} (-1) & e_{
\epsilon_1 + \epsilon_2} (-1) &  e_{ \epsilon_1 + \epsilon_3} (-1)
\\ e_{ \epsilon_1 + \epsilon_2} (-1) & e_{ 2
\epsilon_2} (-1)&  e_{ \epsilon_2 + \epsilon_3} (-1)
\\e_{ \epsilon_1
+ \epsilon_3} (-1) & e_{ \epsilon_2 + \epsilon_3} (-1) & e_{ 2
\epsilon_3} (-1)
\end{array} \right| {\1}=0
 \end{gather*}
in $\widetilde{L}_{C_{\ell} ^{(1)}} (-\Lambda_0)$, so
$\widetilde{L}_{C_{\ell} ^{(1)}} (-\Lambda_0)$ is a certain quotient
of $V_{C_{\ell} ^{(1)}} (-\Lambda_0)$, for $\ell \ge 3$. Thus, any
irreducible module for $\widetilde{L}_{C_{\ell} ^{(1)}}
(-\Lambda_0)$ is an irreducible module for $V_{C_{\ell} ^{(1)}}
(-\Lambda_0)$. The classif\/ication of all irreducible modules in the
category $\mathcal{O}$  for $V_{C_{3} ^{(1)}} (-\Lambda_0)$ was
obtained in \cite[Example 4.1]{A1}. We will apply this result to
obtain a classif\/ication of all irreducible ordinary
$\widetilde{L}_{C_{\ell} ^{(1)}} (-\Lambda_0)$-modules. (Recall
that a module is called ordinary  if $L(0)$ acts semisimply with
f\/inite-dimensional weight spaces.)

\begin{proposition} \label{classif-C-0} \sloppy Let $\ell \ge 3$.
The set
\begin{gather} \label{skup}  \{ L_{C_{\ell} ^{(1)}}( (-n-1 ) \Lambda_0 + n
\Lambda_1) \; \vert \; n \in {\Zp} \}  \cup \{ L_{C_{\ell} ^{(1)}}(- 2
\Lambda_0 + \Lambda_2) \}
 \end{gather}
 provides a complete list of irreducible ordinary modules for the vertex operator algebras $V_{C_{\ell} ^{(1)}}
(-\Lambda_0)$ and $\widetilde{L}_{C_{\ell} ^{(1)}} (-\Lambda_0)$.
\end{proposition}

\begin{proof}
We use the well-known method   for classif\/ication of highest weights
of $V_{C_{\ell} ^{(1)}} (-\Lambda_0)$-modules as solutions of
certain polynomial equations arising from the singular vectors (cf.~\cite{A-1994,AM,AP,MP,P2,P1}). The highest weights of ordinary modules are of the form
$- \Lambda_0 + \mu$, where $\mu = \sum\limits_{i=1}^{\ell}h_{i} \epsilon
_{i}$. Clearly, $h_{i} \in \Zp$, for $i=1, \ldots , \ell$. Using
polynomials from \cite[Example~4.1]{A1}, and the fact that
\[
\tfrac{1}{2}f_{\epsilon _{3}- \epsilon _{i}}(0)^2 \Delta_3(-1){\1}=
\left| \begin{array}{ccc} e_{ 2\epsilon_1} (-1) & e_{ \epsilon_1 +
\epsilon_2} (-1) &  e_{ \epsilon_1 + \epsilon_i} (-1)
\\ e_{ \epsilon_1 + \epsilon_2} (-1) & e_{ 2
\epsilon_2} (-1)&  e_{ \epsilon_2 + \epsilon_i} (-1)
\\e_{ \epsilon_1
+ \epsilon_i} (-1) & e_{ \epsilon_2 + \epsilon_i} (-1) & e_{ 2
\epsilon_i} (-1)
\end{array} \right| {\1},
\]
for $i=4, \ldots , \ell$, we obtain that the weights $\mu$ are
annihilated by the polynomials
\begin{gather*}
 p_i( \mu) = (h_1 +1) (h_2 + \tfrac{1}{2}) h_i,  \nonumber \\
 q_i( \mu) = (h_1 +1) ( 4h_i + ( h_2 + h_i) ( h_2 + h_i -1) ),
\nonumber
\\
 r_i( \mu) = 4 h_i ( h_2 +1) + ( h_1 + h_i -1) ( h_2 + h_i + h_2 (
h_1 + h_i)), \nonumber
\end{gather*}
for $i=3, \ldots , \ell$. Then one easily obtains that $h_i=0$, for
$i=3, \ldots , \ell$, and that either $h_2=0$ and $h_1=n$, for $n
\in {\Zp}$ or $h_2=1$ and $h_1=1$. So we have proved that any
irreducible ordinary $V_{C_{\ell} ^{(1)}} (-\Lambda_0)$-module
(resp. $\widetilde{L}_{C_{\ell} ^{(1)}} (-\Lambda_0)$-module) must
belong to the set~(\ref{skup}).

Since $a_1 ^+ (-1/2) ^{n} {\1} \in M_{2 \ell}$ is a singular vector
of highest weight $-(n+1) \Lambda_0 + n \Lambda_1$, and
\[
e_{\epsilon _{1}+ \epsilon _{2}} ^*=\tfrac{1}{2}(a_1 ^+
(-\tfrac{1}{2}) a_{2\ell -1} ^-( -\tfrac{1}{2}){\1} - a_2 ^+
(-\tfrac{1}{2}) a_{2\ell} ^-( -\tfrac{1}{2}) {\1}) \in M_{2 \ell}
\]
is a singular vector of highest weight $-2\Lambda_0 + \Lambda_2$, we
have that every module from the set~(\ref{skup}) is a module for
these vertex operator algebras. The proof is now complete.
\end{proof}

\begin{remark}
In our new paper \cite{AP-2012} we prove that the vertex operator
algebra $\widetilde{L}_{C_{\ell} ^{(1)} } (-\Lambda_0)$ is simple.
Therefore, Proposition~\ref{classif-C-0} also gives the
classif\/ication of irreducible $L_{C_{\ell} ^{(1)}
}(-\Lambda_0)$-modules.
\end{remark}

\section[Conformal embedding of $C_{\ell} ^{(1)}$ into $A_{2 \ell -1} ^{(1)}$ at level $-1$]{Conformal embedding of $\boldsymbol{C_{\ell} ^{(1)}}$ into $\boldsymbol{A_{2 \ell -1} ^{(1)}}$ at level $\boldsymbol{-1}$}

In this section we show that $L_{C_{\ell} ^{(1)}} (-\Lambda_0)$ is a
$\Z _2$-orbifold of vertex operator algebra $L_{A_{2 \ell -1}
^{(1)}} (-\Lambda_0)$, and determine the corresponding
decomposition.

The subalgebra of $M_{2 \ell}$ generated by
\begin{gather*}  e_{\epsilon _{i}- \epsilon _{j}}^{A} = a_i ^+
(-\tfrac{1}{2}) a_{j} ^-( -\tfrac{1}{2})
 {\1}, \qquad f_{\epsilon _{i}- \epsilon _{j}}^{A}=  a_i ^- (-\tfrac{1}{2}) a_{j} ^+( -\tfrac{1}{2})
 {\1}  \qquad   \mbox{for} \
  i,j =1, \dots, 2 \ell, \quad i<j,
\end{gather*}
is a level $-1$ af\/f\/ine vertex operator algebra associated to the
af\/f\/ine Lie algebra ${\hat \g_1}$ of type $A_{2 \ell -1}^{(1)}$. We
denote it by $\widetilde{L}_{A_{2 \ell -1} ^{(1)}} (-\Lambda_0)$.
Clearly, $\widetilde{L}_{C_{\ell} ^{(1)}} (-\Lambda_0)$ is a
subalgebra of $\widetilde{L}_{A_{2 \ell -1} ^{(1)}} (-\Lambda_0)$.
Let ${\hg}$ be the af\/f\/ine Lie algebra of type $C_{\ell} ^{(1)}$.

As before, denote by $ \omega_1$ the Virasoro vector in ${\widetilde
L}_{C_{\ell} ^{(1)}}(- \Lambda_0)$, and by $ \omega_1 ^A$ the
Virasoro vector in $\widetilde{L}_{A_{2 \ell -1} ^{(1)}}
(-\Lambda_0)$.

\begin{proposition} We have
\[ \omega_1 = \omega_1 ^A.\]
\end{proposition}

\begin{proof} Similarly as in Proposition \ref{Vir-decomp-C}, one can show that
\begin{gather*} \omega _{1}^A= \frac{1}{2} \sum_{i=1} ^{2 \ell} \big( a_i ^- (-
\tfrac{3}{2}) a_i ^{+} (-\tfrac{1}{2} ) -a_i ^+ (- \tfrac{3}{2}) a_i
^{-} (-\tfrac{1}{2} ) \big) {\1} + \frac{1}{4 \ell} H(-1) ^{2} {\1},
 \end{gather*}
 which implies the claim of proposition.
\end{proof}

Furthermore, the vector $ e_{\epsilon _{1}+ \epsilon _{2}} ^*$ from
the proof of Proposition \ref{classif-C-0} is a singular vector for
${\hg}$  in $\widetilde{L}_{A_{2 \ell -1} ^{(1)}} (-\Lambda_0)$
which generates $\widetilde{L}_{C_{ \ell} ^{(1)}}
(-\Lambda_0)$-module $\widetilde{L}_{C_{\ell} ^{(1)}} (-2\Lambda_0
+ \Lambda_2)$, whose top component is irreducible $\frak g$-module
$V_{C_{\ell}}(\omega _2)$. Clearly, $\theta$ acts as $1$ on
$\widetilde{L}_{C_{\ell} ^{(1)}} (-\Lambda_0)$ and as $-1$ on
$\widetilde{L}_{C_{\ell} ^{(1)}} (-2\Lambda_0 + \Lambda_2)$.

\begin{lemma} \label{lem-ext}
Let $u,v \in \widetilde{L}_{C_{\ell} ^{(1)}} (-2\Lambda_0 +
\Lambda_2)$. Then $u_n v \in \widetilde{L}_{C_{\ell} ^{(1)}}
(-\Lambda_0)$, for any $n \in \Z$.
\end{lemma}
\begin{proof} It suf\/f\/ices to prove the lemma for $u$ and $v$ from top component $R(0)$  of $\widetilde{L}_{C_{\ell} ^{(1)}} (-2\Lambda_0 +
\Lambda_2)$. Then the statement will follow from the associator
formulae.

First we notice that \[ u _0 v \in \widetilde{L}_{C_{\ell} ^{(1)}}
(-\Lambda_0),\]
(since  vectors of conformal weight $1$ with bracket $[u,v] =u_0 v$
span   Lie algebra ${\g}_1$ of type $A_{2 \ell -1}$, and f\/ixed point
subalgebra ${\g}$ is a Lie algebra of type $C_{\ell}$).

Assume now that    \[ u_{n_0} v \notin \widetilde{L}_{C_{\ell}
^{(1)}} (-\Lambda_0)\]
for certain $u, v \in R(0)$ and $n_0 \in {\Z}$. Take maximal $n_0$
with this property.

 Then  $u_{n_0}  v$ has nontrivial component  in some highest weight $\widetilde{L}_{C_{\ell} ^{(1)}}
(-\Lambda_0)$-module $W$  of highest weight $-\Lambda_0 + \mu$, and
therefore there is a nontrivial intertwining operator of type
\[ {W \choose \widetilde{L}_{C_{\ell} ^{(1)}} (-2\Lambda_0 +
\Lambda_2) \quad \widetilde{L}_{C_{\ell} ^{(1)}} (-2\Lambda_0 +
\Lambda_2) }. \]
One can associate to this intertwining operator, a non-trivial
$\g$-homomorphism \[ f: \  V_{C_{\ell}}(\omega _2) \otimes
V_{C_{\ell}}(\omega _2)  \rightarrow V_{C_{\ell}}(\mu).\] In
particular, $V_{C_{\ell}}(\mu)$ must appear in the decomposition of
tensor product $V_{C_{\ell}}(\omega _2) \otimes V_{C_{\ell}}(\omega
_2)$.

First consider the case $\ell =2$. Using the decomposition of tensor
product $V_{C_{2}}(\omega _2) \otimes V_{C_{2}}(\omega _2)$ from
(\ref{tens-pr-decomp}) and the fact that the lowest conformal
weights of modules of highest weights $-\Lambda_0 + 2 \omega _2$ and
$-\Lambda_0 + 2 \omega _1$ are $\frac{5}{2}$ and $\frac{3}{2}$,
respectively, we conclude that these modules cannot appear inside
$\widetilde{L}_{A_{2 \ell -1} ^{(1)}} (-\Lambda_0)$. Thus, $u_{n_0}
v \in \widetilde{L}_{C_{\ell} ^{(1)}} (-\Lambda_0)$.

Now, let $\ell \geq 3$. Relation (\ref{tens-pr-decomp}) implies that
the only cases of weights $-\Lambda_0 + \mu$ (aside from $\mu =0$)
from Proposition \ref{classif-C-0} such that $\mu$ appears in the
decomposition of $V_{C_{\ell}}(\omega _2) \otimes
V_{C_{\ell}}(\omega _2)$ are when $\mu = 2 \omega _1$ or $\mu =
\omega _2$. Since  $u_{n_0} v$ is $\theta$-invariant, it does not
contain component inside $\widetilde{L}_{C_{\ell} ^{(1)}}
(-2\Lambda_0 + \Lambda_2)$. On the other hand, the lowest conformal
weight of module of highest weight $-\Lambda_0 + 2 \omega _1$ is
$\frac{\ell + 1}{\ell}$, which is not an integer. Thus, this module
cannot appear inside $\widetilde{L}_{A_{2 \ell -1} ^{(1)}}
(-\Lambda_0)$. Thus, $u_{n_0} v \in \widetilde{L}_{C_{\ell} ^{(1)}}
(-\Lambda_0)$.
\end{proof}

\begin{theorem} \label{thm-ext} We have:
\[\widetilde{L}_{A_{2 \ell -1} ^{(1)}} (-\Lambda_0) =
\widetilde{L}_{C_{\ell} ^{(1)}} (-\Lambda_0) \oplus
\widetilde{L}_{C_{\ell} ^{(1)}} (-2\Lambda_0 + \Lambda_2).\] In
particular, \begin{gather*}  \widetilde{L}_{C_{\ell} ^{(1)}}
(-\Lambda_0)= \widetilde{L}_{A_{2 \ell -1} ^{(1)}} (-\Lambda_0)
^{0}, \qquad
  \widetilde{L}_{C_{\ell} ^{(1)}} (-2\Lambda_0 + \Lambda_2)=
\widetilde{L}_{A_{2 \ell -1} ^{(1)}} (-\Lambda_0) ^{1}.
\end{gather*}
\end{theorem}
\begin{proof}\sloppy
Lemma~\ref{lem-ext} shows that $\widetilde{L}_{C_{\ell} ^{(1)}}
(-\Lambda_0) \oplus \widetilde{L}_{C_{\ell} ^{(1)}} (-2\Lambda_0 +
\Lambda_2)$ is a vertex subalgebra of $\widetilde{L}_{A_{2 \ell -1}
^{(1)}} (-\Lambda_0)$. But this subalgebra clearly contains
generators of $\widetilde{L}_{A_{2 \ell -1} ^{(1)}} (-\Lambda_0)$,
which implies the claim of theo\-rem.
\end{proof}

{\sloppy The classif\/ication of irreducible ordinary $\widetilde{L}_{A_{2 \ell
-1} ^{(1)}} (-\Lambda_0)$-modules follows from the results from~\cite{AP} and similar arguments as in Section~\ref{sect.class.C}:

}

\begin{proposition}\sloppy
The set
\begin{gather*}
\big\{ L_{A_{2\ell -1} ^{(1)}}( (-n-1 ) \Lambda_0 + n \Lambda_1) \;
\vert \; n \in {\Zp} \big\}  \cup \big\{ L_{A_{2\ell -1} ^{(1)}}( (-n-1 )
\Lambda_0 + n \Lambda_{2 \ell -1}) \; \vert \; n \in {\Zp} \big\}
 \end{gather*}
 provides a complete list of irreducible ordinary modules for the vertex operator algebra
 $\widetilde{L}_{A_{2 \ell -1} ^{(1)}} (-\Lambda_0)$.
\end{proposition}

The following result shows that most irreducible
$\widetilde{L}_{A_{2 \ell -1} ^{(1)}} (-\Lambda_0)$-modules remain
irreducible when we restrict them on $\widetilde{L}_{C_{\ell}
^{(1)}} (-\Lambda_0)$.

\begin{proposition}
Assume that $\ell \geq 3$, $n \in {\N}$. Then we have the following
isomorphisms of $\widetilde{L} _{C_{\ell} ^ {(1)} }
(-\Lambda_0)$-modules:
\begin{gather*}
{L}_{A_{2 \ell -1} ^{(1)}} (-(n+1) \Lambda_0 + n \Lambda_1) \cong
L_{C_{\ell} ^{(1)}} (-(n+1) \Lambda_0 + n \Lambda_1), \\
 {L}_{A_{2 \ell -1} ^{(1)}} (-(n+1) \Lambda_0 + n \Lambda_{2\ell -1}
) \cong L_{C_{\ell} ^{(1)}} (-(n+1) \Lambda_0 + n \Lambda_1).
\end{gather*}
\end{proposition}
\begin{proof} We use the Theorem 6.1 from~\cite{DM}.  The def\/inition of automorphism $\theta$ then implies (in
the notation of \cite{DM}) that \[ \theta \circ {L}_{A_{2 \ell -1}
^{(1)}} (-(n+1) \Lambda_0 + n \Lambda_1) \cong {L}_{A_{2 \ell -1}
^{(1)}} (-(n+1) \Lambda_0 + n \Lambda_{2\ell -1} ), \]
as $\widetilde{L}_{A_{2 \ell -1} ^{(1)}} (-\Lambda_0)$-modules.

Theorem 6.1 from \cite{DM} now implies that ${L}_{A_{2 \ell -1}
^{(1)}} (-(n+1) \Lambda_0 + n \Lambda_1)$ and ${L}_{A_{2 \ell -1}
^{(1)}} (-(n+1) \Lambda_0 + n \Lambda_{2\ell -1} )$ are irreducible
as $\widetilde{L} _{C_{\ell} ^ {(1)} } (-\Lambda_0)$-modules. The
claim of Proposition now follows easily.
\end{proof}

 In \cite{AP-2012} we shall prove that   $\widetilde{L}_{A_{2 \ell -1} ^{(1)}} (-\Lambda_0)$ and  $\widetilde{L}_{C_{\ell} ^{(1)}}
(-\Lambda_0)$ are simple. We don't use this result in the present
paper. But even without these simplicity results we can conclude
that the analogous of Theorem \ref{thm-ext} also holds for simple
vertex operator algebras ${L}_{A_{2 \ell -1} ^{(1)}} (-\Lambda_0)$
and $L_{C_{\ell} ^{(1)}} (-\Lambda_0)$.

\begin{corollary} \label{coro-orbifold} We have:
 \[ L_{A_{2 \ell -1} ^{(1)}} (-\Lambda_0) =  L _{C_{\ell} ^{(1)}}
(-\Lambda_0) \oplus  L _{C_{\ell} ^{(1)}} (-2\Lambda_0 +
\Lambda_2).\]
\end{corollary}

\begin{proof}
First we notice that the automorphism $\theta$ also naturally acts
on a simple vertex operator algebra ${L}_{A_{2 \ell -1} ^{(1)}}
(-\Lambda_0)$.  Theorem \ref{thm-ext} implies that \[ {L}_{A_{2 \ell
-1} ^{(1)}} (-\Lambda_0) = V ^0 \oplus V ^1 \] where $V^0$ (resp.
$V^1$) is a quotient of $\widetilde{L}_{C_{\ell} ^{(1)}}
(-\Lambda_0)$ (resp. $\widetilde{L}_{C_{\ell} ^{(1)}} (-2\Lambda_0 +
\Lambda_2)$). By using the fact that ${\Z}_2$-orbifold components
of simple vertex operator algebra are simple (cf.~\cite{DM}), we get
that
$V ^0 =  L _{C_{\ell} ^{(1)}} (-\Lambda_0)$ and $V ^1 =L_{C_{\ell}
^{(1)}} (-2\Lambda_0 + \Lambda_2)$. The proof follows.
 \end{proof}

\section[Commutant of $L _{A_1 ^{(1)}} (-\Lambda_0) ^{ \otimes \ell }$ in
$L_{C_{\ell} ^{(1)}} (-\Lambda_0)$]{Commutant of $\boldsymbol{L _{A_1 ^{(1)}} (-\Lambda_0) ^{ \otimes \ell }}$ in
$\boldsymbol{L_{C_{\ell} ^{(1)}} (-\Lambda_0)}$}

We shall now  study vertex operator algebra $L_{C_{\ell} ^{(1)}}
(-\Lambda_0) $. Corollary \ref{coro-orbifold} implies that
$L_{C_{\ell} ^{(1)}} (-\Lambda_0) $ is a ${\Z}_2$-orbifold of
vertex operator algebra $L_{A_{2 \ell -1} ^{(1)}} (-\Lambda_0)$.
Clearly, $L_{A_{2 \ell -1} ^{(1)}} (-\Lambda_0)$ is a quotient of
$\widetilde{L}_{A_{2 \ell -1} ^{(1)}} (-\Lambda_0)$ modulo the
maximal submodule of $\widetilde{L}_{A_{2 \ell -1} ^{(1)}}
(-\Lambda_0)$ (which is possibly zero).

As before, we denote
\begin{gather*}
H^{(i)}= (a_ i ^+ (-\tfrac{1}{2}) a_ i ^{-} (-\tfrac{1}{2})
+ a_{ 2\ell  +1- i}  ^+ (-\tfrac{1}{2}) a_{2\ell +1- i}  ^ - (-
\tfrac{1}{2}) ){\bf 1} , \quad i = 1, \dots, \ell,
\\
  H = H^{(1)} +  \cdots + H ^{ (\ell)}.
  \end{gather*}
Also, let
\begin{gather*}
 \overline{H}^{(i)} =  H ^{(i)} - H ^{(i+1)}, \quad i =1, \dots, \ell -1, \qquad \mbox{and} \nonumber \\
  \overline{\h}_1 = {\C} \overline{H}^{(1)} + \cdots +
{\C}\overline{H} ^{( \ell-1)}.
\end{gather*}
Clearly $\langle \overline{H} ^{(i)} , \overline{H} ^{(j)} \rangle
=-4 \delta_{i,j} $.

\begin{theorem}
We have:
\[ {\rm Com} \big(   L _{A_1 ^{(1)}} (-\Lambda_0) ^{\otimes \ell}
,L_{C_{\ell}  ^{(1)}} (-\Lambda_0) \big)    \cong M_ { \overline{\h}_1}
(1) ^+.  \]
\end{theorem}
\begin{proof}
By using same arguments as in the proof of Proposition
\ref{com-bos-higher} we get that:
\[   {\rm Com}  \big(   L _{A_1
^{(1)}} (-\Lambda_0) ^{\otimes \ell}  ,\widetilde{ L}_{A_{2 \ell-1}
^{(1)}} (-\Lambda_0) \big) = M_ { \overline{\h}_1}  (1). \]

 Since generators $\overline{H} ^{(i)}$ don't belong to the maximal submodule of  $\widetilde{ L}_{A_{2 \ell-1}  ^{(1)}} (-\Lambda_0) $
  (which is possibly zero) we conclude that
\begin{gather} \label{prva-tv}
{\rm Com} \big(   L _{A_1 ^{(1)}} (-\Lambda_0)
^{\otimes \ell}  , L_{A_{2 \ell -1 }  ^{(1)}} (-\Lambda_0) \big) \cong
M_ { \overline{\h}_1}  (1)  .
\end{gather}

Let $\theta$ be the automorphism of order two of the vertex operator
algebra $L_{A_{2 \ell -1 }  ^{(1)}} (-\Lambda_0) $ as in Corollary
\ref{coro-orbifold}. Then
\[ \theta ( \overline{H}  ^{(i)} ) = -
\overline{H} ^{(i)}, \qquad i = 1, \dots, \ell -1,
\]
and we get that
\begin{gather} \label{druga-tv} M_ {
\overline{\h}_1}  (1) ^+ \subset L_{A_{2 \ell -1 }  ^{(1)}}
(-\Lambda_0)  ^0 = L_{C_{\ell} ^{(1)}} (-\Lambda_0) .
\end{gather}
  Now proof of the theorem follows from relations (\ref{prva-tv}) and (\ref{druga-tv}).
\end{proof}

\subsection*{Acknowledgements}

The authors gratefully acknowledge partial support by the Ministry
of Science, Education and Sports of the Republic of Croatia, Project
ID 037-0372794-2806.

\pdfbookmark[1]{References}{ref}
\LastPageEnding

\end{document}